\documentclass[12pt,a4paper]{article}
\usepackage{times,exscale,amsfonts,amssymb,latexsym,calrsfs}

\def\R{\mathbb R}
\def\C{\mathbb C}
\def\N{\mathbb N}
\def\Z{\mathbb Z}

\makeatletter
\@addtoreset{equation}{section}

\def\@listI{
    \leftmargin\leftmargini
    \parsep 3pt plus 3pt minus 2pt
    \topsep 3pt plus 3pt minus 2pt
    \itemsep \parsep}
\let\@listi\@listI
\makeatother

\def\pm{{\mathsurround=0pt\mathchoice
{\mathop{\raise-3pt\hbox to 0pt {$\displaystyle-$\hss} 
	\raise2.5pt\hbox{$+$}}\nolimits}
{\mathop{\raise-3pt\hbox to 0pt {$\textstyle-$\hss} 
	\raise2.5pt\hbox{$\textstyle+$}}\nolimits}
{\raise-1.8pt\hbox to 0pt {$\scriptscriptstyle-$\hss} 
	\raise1.6pt\hbox{$\scriptscriptstyle+$}}
{\raise-1.8pt\hbox to 0pt {$\scriptscriptstyle-$\hss} 
	\raise1.6pt\hbox{$\scriptscriptstyle+$}}
}}

\def\mp{{\mathsurround=0pt\mathchoice
{\mathop{\raise-2pt\hbox to 0pt {$\displaystyle+$\hss} 
        \raise3.6pt\hbox{$-$}}\nolimits}
{\mathop{\raise-2pt\hbox to 0pt {$\textstyle+$\hss} 
        \raise3.6pt\hbox{$\textstyle-$}}\nolimits}
{\raise1.8pt\hbox to 0pt {$\scriptscriptstyle-$\hss} 
        \raise-1.2pt\hbox{$\scriptscriptstyle+$}}
{\raise1.8pt\hbox to 0pt {$\scriptscriptstyle-$\hss} 
        \raise-1.2pt\hbox{$\scriptscriptstyle+$}}
}}

\newcommand{\ud}{\mathrm{d}}  
\newcommand{\dg}{L^{2}(G)}
\newcommand{\rlb}{$\mathcal{L}(\lambda)$}
\newcommand{\rrlb}{\mathcal{L}(\lambda)}

\newcommand{\sigll}{\mathfrak{S}(\mathcal{L})}
\newcommand{\sigl}{\mathfrak{S}_{0,1}(\mathcal{L})}
\newcommand{\kkk}{K_0^0(\Gamma)}
\newcommand{\kkkk}{K_{1+\varepsilon}^0(\Gamma)}
\newcommand{\ta}{\tilde a}
\newcommand{\tg}{\tilde g}
\newcommand{\tu}{\tilde u}
\newcommand{\loc}{{\mathrm{loc}}}
\renewcommand{\div}{\mathop{\mathrm{div}}}

\newcommand{\sh}{\mathop{\mathrm{sh}}\nolimits}
\newcommand{\ch}{\mathop{\mathrm{ch}}\nolimits}

\newcommand{\argch}{\mathop{\mathrm{argch}}\nolimits}

\usepackage{theorem}
\newtheorem{defi}{Definition}[section]
\newtheorem{theo}[defi]{Th\'eor\`eme}
\newtheorem{cor}[defi]{Corollaire}
\newtheorem{lem}[defi]{Lemme}

\newcommand{\QED}{\mbox{}\hfill\raisebox{-0.2pt}{\rule{5.6pt}{6pt}\rule{0pt}{0pt}} 
          \smallskip\par}             
\newenvironment{Proof}{\noindent
          \parindent=0pt\abovedisplayskip = 0.8\abovedisplayskip
          \belowdisplayskip=\abovedisplayskip{\sc D\'emonstration. }}{\QED}

\renewcommand{\Re}{\mathop{\rm Re}}
\renewcommand{\Im}{\mathop{\rm Im}}
             
\newcommand{\Ref}[1]{\mbox{\rm (\ref{#1})}}            
\newcommand{\REF}[1]{~\mbox{\rm \ref{#1}}}             
\newcommand{\Norm}[2]{\|#1\|\left.\vphantom{T_{j_0}^0}\!\!\right._{#2}}                    
                    
\newcommand{\Res}{\mathop{\mathrm {R\acute es}}}
\newcommand{\reg}{{\mathrm{reg}}}

\title{Probl\`emes de transmission non coercifs dans des polygones}
\author{Monique Dauge et Benjamin Texier}
\date{}

\begin{document}

\maketitle

\begin{abstract}
Nous d\' emontrons que la th\' eorie de Kondratev (1967) pour les domaines \`a coins s'\' etend au cas d'un problme de transmission scalaire entre deux mat\' eriaux de forme polygonale si le contraste entre ces deux mat\' eriaux est n\' egatif et diff\' erent de $-1$. 

\smallskip
{\bf Note.} Ce texte a \'et\'e r\'edig\'e \`a l'issue du stage de Benjamin {\sc Textier} pour sa premi\`ere ann\'ee \`a l'ENS Lyon. Il est paru comme pr\'epublication de l'IRMAR (Universit\'e de Rennes 1) en d\'ecembre 1997 avec le num\'ero {\bf 97-27}. Ce travail ayant maintenant des applications en relation avec les m\'eta-mat\'eriaux, les deux auteurs ont d\'ecid\'e de le d\'eposer sur des serveurs d'archives ouvertes (HAL et arXiv) en f\'evrier 2011.
\end{abstract}


\section{Pr\'esentation du probl\`eme}

Un probl\`eme de transmission scalaire d'ordre $2$ pos\'e sur un domaine plan born\'e
$\Omega$ peut s'exprimer comme probl\`eme variationnel associ\'e \`a une forme
bilin\'eaire $b$ de la forme
$$
   \forall\,u,\,v\in  H^1(\Omega), \quad b(u,v) = \int_{\Omega}\,a\,\nabla u\nabla v,
$$
o\`u $a$ est une fonction \`a valeurs matrices $2\times2$ d\'efinie sur
$\overline\Omega$ et r\'eguli\`ere par morceaux sur une partition finie $(\Omega_j)$
de $\Omega$. Les solutions de tels probl\`emes satisfont des conditions de saut sur
leurs d\'eriv\'ees normales aux interfaces entre les sous-domaines $\Omega_j$ o\`u
$a$ est discontinue. De tels probl\`emes sont li\'es aux
\'equations de Maxwell dans des mat\'eriaux h\'et\'erog\`enes.

Si la fonction $a$ satisfait une condition de positivit\'e uniforme du type 
$$
   \exists \rho_0,\rho_1>0,\quad \forall x\in\overline\Omega,\quad\quad
   \rho_0|\xi|^2 \le  (a(x)\xi \cdot \bar\xi) \le \rho_1|\xi|^2
$$
la forme $b$ est coercive sur $H^1(\Omega)$, et le probl\`eme de transmission peut
\^etre qualifi\'e d'{\it elliptique}.

Par contre se pr\'esente aussi la situation, en liaison avec l'application aux
\'equations de Maxwell pr\'e-cit\'ee (intervention de mat\'eriau supra-conducteur),
o\`u $a$ est positif dans certains sous-domaines et n\'egatif dans les autres. Ici se
pose alors la question de l'ellipticit\'e du probl\`eme de transmission, en plus du
comportement de la solution au voisinage des \'eventuelles irr\'egularit\'es de la
fronti\`ere des $\Omega_j$.

Ici, nous allons, dans un cadre tr\`es simple, \'etudier ces questions, \`a
savoir l'ellipticit\'e la long d'une interface r\'eguli\`ere ainsi que le
comportement au voisinage d'un coin entre deux sous-domaines.
Ainsi, nous consid\'erons comme domaine un rectangle $\Omega$ qui se
d\'ecompose en $2$ sous-domaines $\Omega_{-}$ et $\Omega_{+}$, o\`u
$\overline\Omega_-$ est un rectangle contenu dans $\Omega$, et $\Omega_+$ son
compl\'ementaire.
  
On s'int\'eresse  \`a la r\'esolution du probl\`eme de transmission avec
condition de Dirichlet ext\'erieure\,: 
\begin{equation}
   \label{pb}
   \forall\,v\in  H_0^1(\Omega),\quad \int_{\Omega}\,a\,\nabla u\nabla v - f\,v =0,
\end{equation}
 o\`u $a$ est une fonction constante par morceaux: 
$$  a = \left\{ \begin{array}{ll}
   a_{+} \,\mathrm {Id} & \mbox{ dans } \Omega_{+} \\
   a_{-} \,\mathrm {Id} & \mbox{ dans } \Omega_{-} , \\
              \end{array} \right.
$$ 
et o\`u le second membre $f$ est suppos\'e appartenir \`a $L^{2}(\Omega)$.

 L'\'etude qui suit traite le cas $a_{-}\,a_+<0$, le cas $a_{-}\,a_+>0$ \'etant bien
connu, voir {\sc Nicaise} \cite{NicaiseBook}.  Etant donn\'ee $u\in H_0^{1}(\Omega)$
solution de \Ref{pb}, on se propose d'en \'etudier la structure, la difficult\'e
venant de $a_{-}\,a_+<0$ et de la pr\'esence des coins.

 On notera $u_+$ pour $u|_{\Omega_+}$ et $u_-$ pour $u|_{\Omega_-}$.  D'apr\`es
(\ref{pb}):
$$
   \forall\,v \in \mathcal{D}(\Omega_+),\quad - \int_{\Omega_+}\,a_+\,\Delta u_+\,v 
   = \int_{\Omega_+} f_+\,v.
$$
 Donc:
  $$ - a_+ \Delta u_+ = f_+,$$
et, de m\^eme: $ - a_- \Delta u_- = f_-.$
Alors, si on int\`egre (\ref{pb}) par parties, on obtient, pour tout
$v\in \mathcal{D}(\Omega)$ la condition suivante sur l'interface $\partial\Omega_+
\cap \partial\Omega_-$ (qui est \'egale \`a $\partial\Omega_-$) entre $\Omega_+$ et
$\Omega_-$:
\begin{equation}
  \label{formelle}
  \int_{\partial\Omega_-}
  ( a_+ \partial_{n_+}u_+  + a_- \partial_{n_-}u_-)\,v = 0.
\end{equation}
La solution $u$ n'\'etant a priori que dans $H^1(\Omega)$, (\ref{formelle}) a
seulement un sens comme dua\-lit\'e entre $H^{-1/2}(\partial\Omega_-)$ et
$H^{1/2}(\partial\Omega_-)$.  L'\'etude qui suit montre qu'en r\'ealit\'e, $u$ est
$H^2$, sauf peut-\^etre au voisinage des coins, ce qui donne un sens usuel de trace
aux d\'eriv\'ees normales.

Apr\`es le choix d'une normale ($\partial_n =\partial_{n_+}$, par exemple), 
on obtient la condition dite de transmission:
$$ a_+ \,\partial_{n}u_+  \,-\, a_-\,\partial_{n}u_- = 0,\quad 
\mbox{sur $\partial\Omega_-$}.
$$

Il r\'esulte de l'ellipticit\'e du Laplacien que $u_+$ est dans $H^2_\loc(\Omega_+)$
et $u_-$ dans $H^2_\loc(\Omega_-).$ La condition de Dirichlet \'etant elliptique, on a
aussi la r\'egularit\'e $H^2$ pour $u$ au voisinage de tout point du bord externe
$\partial\Omega$ qui ne rencontre pas l'un des sommet de $\Omega$. De plus, comme les
angles de $\Omega$ sont $\pi/2$, le principe de r\'eflexion (et aussi la th\'eorie
g\'en\'erale du Laplacien dans un polygone, voir {\sc Grisvard} \cite{GrisvardBook})
montre que $u$ est aussi $H^2$ au voisinage des coins de $\Omega$. Il reste donc \`a
\'etudier $u$ au voisinage des points de l'interface $\partial\Omega_-$.

 Dans la section 2, on va montrer que si $a_+ a_- \neq -1$, $u_+$ et $u_-$ ont la
r\'egularit\'e $H^2$ au voisinage de tout point r\'egulier de l'interface. Dans la
section 3, on va montrer que, sous la m\^eme condition $a_+ a_- \neq -1$, le
comportement de $u$ au voisinage des coins de $\Omega_-$ (interface \`a coin)
rel\`eve de la th\'eorie g\'en\'erale des probl\`emes elliptiques dans les ouverts
\`a coins: selon la valeur du rapport $\mu := a_+/a_-$, on obtiendra, soit la
r\'egularit\'e $H^2$ de $u_+$ et $u_-$, soit une d\'ecomposition en parties
r\'eguli\`ere et singuli\`ere. Nous concluons en section 4 par quelques persepectives.


\section{Etude au voisinage d'une interface} 

On se place maintenant au voisinage d'un point de
$\partial\Omega_-$ qui ne rencontre aucun coin. Dans une premi\`ere partie, on
tronque $u$ pour pouvoir se placer dans $\R^2$ tout entier. Puis, par
transformation de Fourier partielle et un changement de variable, on aboutit \`a
l'\'etude d'un probl\`eme \`a une dimension. On d\'eduit de son \'etude une
propri\'et\'e de r\'egularit\'e locale pour $u$.

\subsection{Localisation}
\label{2SS1}
 
Soit $(x_0,y_0)$ un point de l'interface: on suppose que l'interface est localement
contenue dans la droite d'\'equation $y=y_0$. Soit $\chi_0 \in \mathcal{D}(\R)$,
valant $1$ au voisinage de $0$ et \`a support suffisamment petit et soit
$\chi\in\mathcal{D}(\R^2)$ d\'efini par $\chi(x,y) = \chi_0(x - x_0)\,\chi_0( y -
y_0)$. Alors, pour $v \in H^1(\Omega)$:
$$
   \int_{\Omega}a\nabla (\chi u)\nabla v = 
   \int_{\Omega} a \bigl( \nabla u\nabla(\chi v) + 
   u \nabla \chi\nabla v - v\nabla u\nabla\chi\bigr). 
$$ 
En int\'egrant par parties:
$$
   \int_{\Omega_+}a_+ u \nabla\chi \nabla v = 
   - \int_{\Omega_+} \div(a_+ u\nabla\chi)\,v 
   + \int_{\partial\Omega_+} a_+ uv\,\partial_n\chi,
$$ 
et de m\^eme sur $\Omega_-$. Par choix de $\chi$, $\partial_n\chi = 0$. On a donc:
$$
   \int_{\Omega}\,a\,\nabla (\chi u)\,\nabla v =  
   \int_{\Omega} \bigl(\chi f +  a \nabla u\nabla\chi 
   +  \div(a\, u\nabla\chi) \bigr)\,v.
$$ 
Or, $\div(a\, u\nabla\chi) = \partial_x(a \,u \,\partial_x \chi)
+ \partial_y(a\, u\, \partial_y \chi) = a\bigl(\partial_x(u \,\partial_x \chi) +  
\partial_y( u \,\partial_y \chi)\bigr)$, car $\partial_y\chi=0$ \`a l'interface.
Par suite, comme $u\in H^1(\Omega)$, la fonction $g$ d\'efinie par $g = \chi f +  a
\nabla u\nabla\chi +  \div(a \,u\nabla\chi)$ est dans $L^2(\Omega)$. Ainsi:
$$
   \forall\,v\in H_0^1(\Omega),\quad \int_{\Omega}a\nabla (\chi u)\nabla v =
   \int_{\Omega} g \,v .
$$ 
Donc, si on note $\tilde u$ le prolongement de $\chi u$ par $0$ \`a tout $\R^2$ et de
m\^eme $\tilde g$ le prolongement de $g$ par $0$, et en posant
$(x_0,y_0)=(0,0)$, 
et  $\tilde a = \left\{ \begin{array}{cc}  a_+ & \mbox{pour $y>0$} \\
                                   a_- & \mbox{pour $y<0$} 
                      \end{array} \right.$:
\begin{equation}
\label{pb2}
   \forall\,v\in H^1(\R^2),\quad 
   \int_{\R^2}\,\tilde a\,\nabla \tu\,\nabla v =
   \int_{\R^2} \tg \,v .
\end{equation}

\subsection{R\'egularit\'e}
\label{2SS2}
On part maintenant de l'\'equation \Ref{pb2}, et on renote $\tu$ par $u$ et $\tg$ par
$g$ dans cette sous-section\REF{2SS2}. En proc\'edant comme dans la section 1, on voit
que $u$ v\'erifie les \'equations:
$$\left \{ \begin{array}{ll}
               -a_+\,\Delta u_+ = g_+ & \mbox{ pour } y>0, \\[1mm]
               -a_-\,\Delta u_- = g_- & \mbox{ pour } y<0. 
            \end{array}  \right.   
$$
Appliquons la transformation de Fourier partielle par rapport \`a $x$: 
$$
   a_+\, (\xi^2 - \partial_y^2)\,\hat{u}_+(\xi,y) = \hat{g}_+(\xi,y),
   \quad\forall\,\xi\in \R,y\in \R_+,
$$ 
de m\^eme pour $y<0$. En effectuant le changement de variable $\tilde y=|\xi|\, y$
(et o\`u $\tilde y$ est imm\'ediatement renot\'e $y$), on a\,:
$$
   a_\pm\,(1 - \partial_y^2)\,\hat{u}_\pm(\xi,\frac{y}{|\xi|}) =
   \frac{1}{\xi^2}\,\hat{g}_\pm(\xi,\frac{y}{|\xi|}).
$$
Le changement de variable pr\'ec\'edent permet ainsi de fixer $\xi=1$. On consid\`ere
maintenant un probl\`eme \`a une dimension, la variable \'etant $y$.

Notons $h(y) = \frac{1}{\xi^2}\hat{g}(\xi,\frac{y}{|\xi|})$, et 
$E= H^1(\R)\cap \bigl( H^2(\R_-)\times  H^2(\R_+)\bigr)$ i.e.
$$
   E = \big\{ w\in H^1(\R)\;|\quad 
   w_-\in H^2(\R_-)\quad\mbox{et}\quad w_+\in  H^2(\R_+)\bigr\}.
$$
On s'int\'eresse au probl\`eme\,:
\begin{equation}
\label{eq1D}
   \mbox{ Trouver $w\in E$ tel que } \forall\, v \in H^1(\R),\quad
   \int_{\R}\,\ta\,(w\, v + \nabla w\,\nabla v) - h\,v\, = 0.
\end{equation}
 
Supposons conna\^{\i}tre $w\in E$ v\'erifiant (\ref{eq1D}). Alors, en choisissant $v$ 
dans $\mathcal{D}(\R_+)$ et en int\'egrant par parties, on obtient comme dans la
section 1:
$$
   a_+\,(w_+ - w_+'') = h_+ \quad\mbox{et}\quad a_-\,(w_- - w_-'') = h_-.
$$
En reprenant (\ref{eq1D}), on obtient la condition de transmission:
$$
   a_+\, w_+'(0) \, - \, a_-\, w_-'(0)\, =\, 0\,.
$$
Consid\'erons alors la solution 
$w_0 = (w_{0,-},w_{0,+})\in H_0^1(\R_-)\times H_0^1(\R_+)$ du probl\`eme de
Dirichlet: 
$$
   a_-(w_{0,-} - w_{0,-}'') = h_-\,,\quad a_+(w_{0,+} - w_{0,+}'') = h_+.
$$
On sait que $w_{0,-} \in H^2(\R_-)$ et $w_{0,+} \in H^2(\R_+)$.
Si on pose $w_1 = w - w_0$, alors n\'ecessairement:
$$ \left \{ \begin{array}{l}
    w_{1,-} - w_{1,-}'' = 0, \quad w_{1,+} - w_{1,+}'' = 0, \\[1mm]
    w_{1,-}(0) = w_{1,+}(0),\\[1mm]
    a_+\,w_{1,+}'(0)  - a_-\, w_{1,-}'(0)  =   
    a_-\, w_{0,-}'(0) - a_+\,w_{0,+}'(0). 
            \end{array}   \right. 
$$
En posant $\alpha = a_-\, w_{0,-}'(0) - a_+\,w_{0,+}'(0)$, on a n\'ecessairement,
si $a_+ + a_- \ne 0\,:$ 
\begin{equation}
\label{ex}
   w_{1,\pm} = \frac{\alpha}{a_+ + a_-}\,e^{\mp y}.
\end{equation}
D'o\`u l'unicit\'e de la solution de \Ref{eq1D}.
R\'eciproquement, en posant (\ref{ex}) et $w = w_0 + w_1$, $w$ est solution de 
(\ref{eq1D}), d'o\`u l'existence. 

Alors, notant $\tilde E = H^{1}(\R)\cap\,\bigl\{ v \in H^2(\R_-)\times 
H^2(\R_+)\,,\,a_+\,v_+'(0) - a_-\,v_-'(0) = 0 \bigr\}$, si $a_+ + a_- \ne 0$:
$$
  \begin{array}{ccc}
              \tilde E & 
          \longrightarrow &
               L^2(\R_-)\times  L^2(\R_+) \\ 
                v &
          \longrightarrow &
                (v_- - v_-''\,,\,v_+ - v_+'')
  \end{array} 
$$
est lin\'eaire, continue, bijective d'apr\`es ce qui pr\'ec\`ede, donc ouverte. 
D'apr\`es le th\'eor\`eme de Banach: 
\begin{equation}
\label{Ban}
   \exists\,C>0,\quad \Norm{ w_{\pm}}{H^2(\R_{\pm})}\,\leqslant\, C\, 
   \bigl( \Norm{ h_-}{L^2(\R_-)} + \Norm{ h_+}{L^2(\R_+)} \bigr).
\end{equation} 

Avec (\ref{Ban}), il suffit maintenant de remonter les calculs pr\'ec\'edents pour 
obtenir le r\'esul\-tat de r\'egularit\'e annonc\'e. En se pla\c{c}ant par exemple
au-dessus de l'interface:
$$
   \Norm{w_+}{H^2(\R_+)}^2 \leqslant C \,\Norm{h}{L^2(\R)}^2. 
$$
On applique ceci pour $w(y) = \hat{u}(\xi,\frac{y}{|\xi|})$ et $h(y) =
\frac{1}{\xi^2}\hat{g}(\xi,\frac{y}{|\xi|})$. Or, 
$$
   \Norm{\partial_y^j w_+}{L^2(\R_+)}^2 = |\xi|^{1-2j} \int_{\R_+}
   |\partial_y^j \hat{u}|^2\, \ud y \quad\mbox{pour}\quad j=0,1,2.
$$ 
Donc  
\begin{equation}
\label{2E3}
   \Norm{\xi^2\hat{u}}{L^2(\R_+)} + \Norm{\xi \partial_y \hat{u}}{L^2(\R_+)} + 
   \Norm{\partial_y^2\hat{u}}{L^2(\R_+)} \leqslant C\,\Norm{\hat{g}}{L^2(\R)} .
\end{equation}
Or la caract\'erisation de $H^2(\R^2_+)$ par transformation de Fourier partielle est
\begin{equation}
\label{2E1}
   \Norm{v}{H^2(\R^2_+)}^2 \simeq 
   \int_{\R_+} \Norm{v(\xi,\cdot)}{H^2(\R_+,\,|\xi|)}^2 \, \ud\xi
\end{equation}
o\`u la norme \`a param\`etre $\Norm{\,\cdot\,}{H^2(\R_+,\,\rho)}$ est d\'efinie pour
$\rho\geqslant0$ et $z\in H^2(\R_+)$:
\begin{equation}
\label{2E2}
   \Norm{z}{H^2(\R_+,\,\rho)}^2 = \Norm{(1+\rho)^2 z}{L^2(\R_+)}^2 +  
   \Norm{(1+\rho)\, \partial z}{L^2(\R_+)}^2 + \Norm{\partial^2 z}{L^2(\R_+)}^2.
\end{equation}
L'espace $H^1(\R^2_+)$ a bien s\^ur une caract\'erisation analogue.
Pour $\rho$ grand:
$$
   \Norm{z}{H^2(\R_+,\,\rho)}^2 \simeq \Norm{\rho^2 z}{L^2(\R_+)}^2 +  
   \Norm{\rho\, \partial z}{L^2(\R_+)}^2 + \Norm{\partial^2 z}{L^2(\R_+)}^2.
$$
Pour les $\rho$ petits, il reste un terme $H^1$ \`a droite. Donc \Ref{2E3} donne:
$$
   \Norm{\hat{u}(\xi,\cdot)}{H^2(\R_+,\,|\xi|)}^2 \leqslant C\,
   \bigl( \Norm{\hat{g}(\xi,\cdot)}{L^2(\R)}^2 +  
   \Norm{\hat{u}(\xi,\cdot)}{H^1(\R_+,\,|\xi|)}^2 \bigr) .
$$ 
Ainsi, d\'eduit-on de \Ref{2E1} que:
\begin{equation}
\label{2E4}
   \Norm{u}{H^2(\R^2_+)}^2 \leqslant C\,\bigl( \Norm{g}{L^2(\R^2)}^2 + 
   \Norm{u}{H^1(\R^2_+)}^2 \bigr).
\end{equation} 

Revenant \`a la situation ant\'erieure de $u$ solution du probl\`eme \Ref{pb}, et
combinant \Ref{2E4}, valable pour $\chi u$, avec les estimations elliptiques
classiques au voisinage des points int\'e\-rieurs de $\Omega_\pm$ et des points
r\'eguliers de $\partial\Omega$, et les estimations de r\'egularit\'e au voisinage
des coins externes de $\partial\Omega$, nous obtenons:

\begin{theo}
  \label{final3}
On suppose que $a_+ a_- \ne -1$.
Soit $u\in H_0^1(\Omega)$, solution du probl\`eme \Ref{pb}.
Soit $\Omega_1$ un ouvert de $\R^2$ tel que $\overline\Omega_1$ ne rencontre pas les
coins de transmission (les coins de $\Omega_-$), et soit $\Omega_2$ tel que 
$\overline \Omega_1 \subset \Omega_2$. On note $\Omega_{1,\pm} = \Omega_\pm \cap
\Omega_1$. Alors $u = (u_-,u_+)\in H^2(\Omega_{1,-})\times H^2(\Omega_{1,+})$, avec
les estimations:
$$
   \Norm{u_\pm}{H^2(\Omega_{1,\pm})}\, \leqslant\, C\,
   \bigl( \Norm{f}{L^2(\Omega_2)} + \Norm{u}{H^1(\Omega_2)} \bigr) .
$$   
\end{theo}


\section{Etude au voisinage d'un coin int\'erieur}

Dans une premi\`ere partie, on tronque $u$ \`a nouveau, afin de se placer au voisinage
d'un coin de transmission, o\`u $\Omega_-$ et $\Omega_+$ co\"{\i}ncident
respectivement avec un voisinage de $0$ des secteurs $\Gamma_-$ et $\Gamma_+$ qui
forment une d\'ecomposition de $\Gamma=\R^2$: en coordonn\'ees polaires
$(r,\theta)$,
$$
   \Gamma_- = \bigl\{(x,y),\,r\in\R_+,\,\theta\in (0,\pi/2)\bigr\}, 
   \quad \Gamma_+ = \bigl\{(x,y),\,r\in\R_+,\,\theta\in (\pi/2,2\pi) \bigr\}.
$$

Puis, avec $h = r^2g$ o\`u $g\in L^2(\Gamma)$ est un nouveau second membre apr\`es
localisation, l'\'etude du symbole Mellin $\mathcal{L}(\lambda)$ de notre probl\`eme
de transmission au voisinage de $0$ va fournir un $u_0$ de
r\'egularit\'e $H^2$ dans un voisinage de $0$, tel que:  
$$
   \mathcal{L}(\lambda)\,\hat{u}_0(\lambda) = \hat{h}(\lambda) 
$$
pour tout $\lambda$ sur une certaine droite $\Re\lambda=$ constante. Une
extension de la formule des r\'esidus permet alors comparer $u_0$ et $u$. 

Dans toute la section 3, on note $z_- = z|_{\Gamma_-}$ et $z_+ = z|_{\Gamma_+}$,
pour toute fonction $z$ d\'efinie sur $\R^2$, et aussi $G_- = (0,\pi/2)$, $G_+ =
(\pi/2,2\pi)$ et $G$ le cercle unit\'e $\R/2\pi\Z$.

\subsection{Localisation}
 Soit $r \mapsto \chi(r),\,\chi \in \mathcal{D}(\overline\R_+)$, \`a support compact,
et qui vaut $1$ au voisinage de $0$. Alors, $\partial_\theta\chi = 0$. Avec un tel
$\chi$, on peut reprendre les calculs de la sous-section\REF{2SS1}, pour obtenir:
$$
   \int_{\Omega}\,a\,\nabla (\chi u)\,\nabla v =
   \int_{\Omega} \bigl(\chi f +  a \nabla u\nabla\chi + \div(a\, u\nabla\chi) 
   \bigr)\,v.
$$ 
Or, $ \div(a\,u\nabla\chi) = \partial_r(a \,u \,\partial_r \chi)$ par choix de
$\chi$,  donc
$g = \chi f +  a \nabla u\nabla\chi +  \div(a \,u\nabla\chi)$ est dans $L^2(\Omega)$.
Apr\`es prolongement par z\'ero\,:
\begin{equation}
\label{pb3}
   \forall\,v\in H^1(\R^2),\quad \int_{\R^2}\,\tilde a\,\nabla \tilde u\,\nabla v
   = \int_{\R^2} \tilde g \,v \,,
\end{equation}
o\`u $\ta$ est \'egal \`a $a_+$ sur $\Gamma_+$ et \`a $a_-$ sur $\Gamma_-$.

\subsection{Transformation de Mellin et espaces \`a poids}
La transformation de Mellin est une transformation de Fourier--Laplace (i.e.\ \`a
argument complexe) radiale. Elle caract\'erise particuli\`erement bien les espaces
\`a poids suivants dont l'introduction est classique dans la th\'eorie des
probl\`emes \`a coins, voir \textsc{Kondrat'ev} \cite{Kondratev}. Dans cette
sous-section, $\Gamma$ d\'esigne un secteur plan:
$$
   \Gamma = \bigl\{(x,y),\,r\in\R_+,\,\theta\in G\bigr\}.
$$
  
\begin{defi}
   \label{esp.poids}
Soit $s\in \N$ et $\gamma\in \R$, 
$$
   K_{\gamma}^{s}(\Gamma) = \bigl\{ v\in L_\loc^{2}(\Gamma) ,
   \quad r^{|\alpha|-s+\gamma}\partial_{x}^{\alpha}v \in L^{2}(\Gamma) , \;
   \forall\,\alpha \in\N^{2} , |\alpha|\,\leqslant\,s \bigr\} 
$$
avec la norme  $\Norm{v}{K_{\gamma}^{s}(\Gamma)}^{2} = \sum_{| \alpha |
\le s}\Norm{r^{|\alpha|-s+\gamma}\partial^{\alpha}v}{L^{2}(\Gamma)}^2$. 
\end{defi}

     Pour $v \in \mathcal{D}(\R_{+})$, la transformation de Mellin de $v$
est not\'ee $\hat{v}$ et est d\'efinie pour tout $\lambda \in \C$ par la
formule
$$
   \hat{v}(\lambda) = \int_{0}^{+\infty}\,r^{-\lambda}\,v(r)\,\frac{\ud r}{r}\, .
$$
Si on pose $t=\log r$, $\breve v(t)=v(r)$, et $\lambda=\xi+ i\eta$, alors 
$\hat v(\lambda)$ est la transform\'ee de Fourier de $t\mapsto e^{-\xi t}\,\breve
v(t)$ calcul\'ee en $\eta$.

       Si $v \in \mathcal{D}(\overline{\Gamma}\setminus\{0\})$, $v$ s'\'ecrit en
coordonn\'ees polaires $v(x)=\tilde v(r,\theta)$, et on d\'efinit alors la
transform\'ee de Mellin de $v$ par la formule:
$$
   \hat{v}(\lambda,\theta) = \int_{0}^{+\infty}\,r^{-\lambda}\, \tilde
   v(r,\theta)\,\frac{\ud r}{r} \,.
$$    
L'\'egalit\'e de Parseval permet de d\'emontrer:

\begin{lem}
  \label{lem1.Mel.poids}   
Soit $\gamma \in\R$ et $\xi=-\gamma-1$. Si $v\in K_{\gamma}^{0}(\Gamma)$, alors la
transformation de Mellin de $v$ est bien d\'efinie pour tout $\lambda = \xi+i\eta$,
$(\eta\in\R)$, et
$$
   (\eta,\theta)\mapsto \hat{v}(\xi+i\eta,\theta) \in L^{2}(\R\times G).
$$
La r\'eciproque est vraie et on a l'\'equivalence: 
$$
   \Norm{v}{K_{\gamma}^{0}(\Gamma)} \simeq
   \bigl(\int_{\R}\,\Norm{\hat{v}(\xi+i\eta)}{L^{2}(G)}^2\,\ud \eta\bigr)^{1/2}. 
$$
\end{lem}

En \'etudiant l'effet du changement de variables $x\mapsto (r,\theta)\mapsto 
(t,\theta)$ sur les op\'erateurs de d\'erivation, on peut \'etendre le lemme
pr\'ec\'edent aux espaces $K_{\gamma}^{s}(\Gamma)$ pour $s \in\N$: 
\begin{lem}
\label{5.4}
Soit $\gamma \in\R$ et $s\in\N$, on pose $\xi=s-\gamma-1$.
\begin{itemize}
  \item
    Si  $v\in K_{\gamma}^{s}(\Gamma)$, alors 
     $\eta\mapsto \hat{v}(\xi+i\eta,\cdot) \in  L^{2}(\R,H^{s}(G))$.\\
    Plus pr\'ecis\'ement, on a l'\'equivalence:
     $$ \Norm{v}{K_{\gamma}^{s}(\Gamma)} \simeq   \bigl(\int_{\R}\,
     \Norm{\hat v(\xi+i\eta)}{H^{s}(G,\,|\eta|)}^2\,\ud \eta\bigr)^{1/2} $$ 
     o\`u pour $\rho>0$,
     $$\Norm{v}{H^{s}(G,\,\rho)}^{2} = \sum_{\beta_1+\beta_2 = s}
       \Norm{(1+\rho)^{\beta_1}\partial^{\beta_2}v}{L^{2}(G)}^{2}.$$
  \item
    R\'eciproquement, si $V(\lambda)$ est d\'efinie pour $\Re\lambda=\xi$ de sorte que
    $\eta \mapsto V(\xi+i\eta,\cdot)$ soit dans $L^{2}(\R,H^{s}(G))$ avec la
    condition
    $$ \int_{\R}\,\Norm{V(\xi+i\eta)}{H^{s}(G,\,|\eta|)}^2\,\ud \eta < \infty $$
    alors, pour tout $\lambda$ de partie r\'eelle \'egale \`a $\xi$, $V(\lambda)$ est 
   la transform\'ee de Mellin d'une fonction $v\in K_{\gamma}^{s}(\Gamma)$, et on a:
    \begin{equation}
  \label{Mel.inv}
    \tilde v(r,\theta) =
    \frac{1}{2\pi}\int_{\R}\,r^{\lambda}V(\lambda,\theta)\,\ud \eta.
    \end{equation}
\end{itemize}
\end{lem}

\subsection{Inversibilit\'e du symbole Mellin}
\label{3S3}

On repart maintenant du probl\`eme \Ref{pb3} o\`u on renote $\tu$ par $u$ et $-\tg$
par $g$. Int\'egrant par parties, on obtient que:
\begin{equation}
\label{pb4}
   \left\{\begin{array}{ll}
   a_\pm \Delta u_\pm = g_\pm & \mbox{dans }\; \Gamma_\pm,\\[1mm]
   a_+ \partial_n u_+ - a_- \partial_n u_- = 0 & \mbox{sur }\; 
   \partial\Gamma_+ \cap \partial\Gamma_-,
\end{array}\right.
\end{equation}
avec $u\in H^1(\Gamma)\cap \bigl(H^2_\loc(\overline\Gamma_-\setminus\{0\}) \times
H^2_\loc(\overline\Gamma_+\setminus\{0\}) \bigr)$, {\it cf} Th\'eor\`eme\REF{final3},
et $g\in L^2(\Gamma)$, tous deux \`a support compact.
 
En appliquant la transformation de Mellin au probl\`eme \Ref{pb4}, on obtient
l'expression de \rlb:
$$
  \begin{array}{cccc}
             \rrlb \; : & \hat E & 
          \longrightarrow &
               L^2(G_-)\times  L^2(G_+) \\ 
               & W &
          \longrightarrow &
                \bigl( a_-( \lambda^2 + \partial_\theta^2)W_-\,,\,
                 a_+( \lambda^2 + \partial_\theta^2)W_+ \bigr)
  \end{array} 
$$
o\`u $\hat E = H^1(G)\cap\,\bigl\{ W \in H^2(G_-)\times 
H^2(G_+)\,,\,a_+\,W_+' - a_-\,W_-' = 0 \mbox{ sur }\partial G_-\cap \partial G_+
\bigr\}$.

On se fixe $\lambda\in\C$. On \'etudie l'inversibilit\'e de $\rrlb$.

\noindent\textbf{Injectivit\'e:}  Soit $W\in \hat E$ tel que $\rrlb W = 0$.
On a donc les conditions de transmission:
$$ \left \{ \begin{array}{rcl}  W_+(2\pi) - W_-(0) &=& 0, \\[1mm]
                                a_+ W_+'(2\pi) - a_-W_-'(0) &=& 0, \\[1mm]
                                W_+(\omega) - W_-(\omega) &=& 0, \\[1mm] 
                                a_+ W_+'(\omega) - a_-W_-'(\omega) &=& 0, 
            \end{array}
                              \right. 
$$
o\`u $\omega=\frac\pi2$, la premi\`ere et la troisi\`eme provenant du fait que $W\in
H^1(G)$. Par ailleurs, $\rrlb W = 0$ donne aussi\,: 
$$
   \forall\,\theta \notin \{0,\frac{\pi}{2},2\pi\},\quad 
   (\lambda^2 + \partial_\theta^2) W = 0,
$$ 
donc il existe $\alpha_-,\alpha_+,\beta_-,\beta_+$ tels que $W_-(\lambda) =
\alpha_-\cos (\lambda\,\theta) +\beta_- \sin (\lambda\,\theta)$ et de m\^eme
$W_+(\lambda) =
\alpha_+\cos (\lambda\,\theta) +\beta_+ \sin (\lambda\,\theta)$.
 En injectant ces solutions dans le syst\`eme pr\'ec\'edent, on obtient:
 $\rrlb$ est injectif $\Longleftrightarrow \mathcal{A}(\lambda) \neq 0$, o\`u, en
posant
$\mu= a_+/a_-$:
$$
   \mathcal{A}(\lambda) = \left|\begin{array}{cccc}
   -1 & 0 & \cos 2\pi\lambda & \sin 2\pi\lambda \\
   0 & -1 & -\mu\,\sin 2\pi\lambda & \mu\, \cos 2\pi\lambda \\
   -\cos\lambda\omega & -\sin\lambda\omega & \cos\lambda\omega & \sin\lambda\omega\\
   \sin\lambda\omega & - \cos\lambda\omega & -\mu\, \sin\lambda\omega & \mu\,
      \cos\lambda\omega 
   \end{array}\right|.
$$
On obtient imm\'ediatement que
$$
   \mathcal{A}(\lambda) = (\mu^2+1) \,\sin\lambda\omega\, \sin\lambda(2\pi-\omega) +
   2\mu\, \bigl(1-\cos\lambda\omega\, \cos\lambda(2\pi-\omega)\bigr),
$$
dont on peut montrer qu'il est \'egal \`a
$$
   \mathcal{A}(\lambda) = (b+c)(b-c),\quad\mbox{o\`u}\quad
   b = (\mu+1)\,\sin\lambda\pi \;\mbox{ et }\;
   c = (\mu-1)\,\sin\lambda(\pi-\omega).
$$

\noindent{\textbf {Surjectivit\'e:}} On aboutit au m\^eme d\'eterminant
d'ordre deux. Donc $\rrlb$ est injectif si et seulement si $\rrlb$ est
surjectif.  On a ainsi obtenu
\medskip

\begin{lem}
\label{3L2}
Si $\lambda\ne0$, $\rrlb$ est bijectif si et seulement si
\begin{equation}
\label{3E1}
   (\mu+1)\,\sin\lambda\pi \;\ne\; \pm\, (\mu-1)\,\sin\lambda(\pi-\omega),
   \quad\mbox{avec}\quad \mu = a_+/a_-.
\end{equation}
\end{lem}

Ici $\omega=\pi/2$ et cette valeur particuli\`ere permet de montrer la propri\'et\'e
suivante, qui sera utile pour la d\'ecomposition de $u$:
\medskip

\begin{lem}
\label{3L3}
Si $\mu$ est r\'eel et $\mu\ne1$, pour tout $\lambda,\, \Re\lambda=1,\,\rrlb $ est
inversible.
\end{lem}
\medskip

\begin{Proof}
En effet pour $\lambda=1+i\eta$, la condition \Ref{3E1} s'\'ecrit
$$
  - i (\mu+1)\,\sh\eta\pi \;\ne \;\pm\, (\mu-1)\,\ch(\eta\pi/2),
$$ 
dont la seule racine est $\mu=1$, $\eta=0$ (probl\`eme de transmission
trivial).
\end{Proof}

\subsection{Calcul des r\'esidus}

Pour analyser une solution $u$ du probl\`eme \Ref{pb} au voisinage d'un coin de
transmission, on se place dans le cadre du probl\`eme \Ref{pb4}, auquel on applique
la transformation de Mellin. Pour cela, on se ram\`ene \`a des espaces \`a poids
$K^s_\gamma[\Gamma)$, voir D\'efinition\REF{esp.poids}. On remarque d'abord que,
comme cons\'equence de l'in\'egalit\'e de Hardy, les fonctions dans $H^1(\Gamma)$ \`a
support compact appartiennent \`a $K_\varepsilon^1(\Gamma)$ pour tout
$\varepsilon>0$. 

On se fixe donc $\varepsilon>0$, que l'on pourra choisir aussi petit que l'on veut
et donc, $u\in K_\varepsilon^1(\Gamma)$. Alors, par le Lemme\REF{5.4}, appliqu\'e \`a
$\Gamma_-$ et $\Gamma_+$ on obtient que pour tout
$\lambda,\,\Re\lambda=-\varepsilon$, $\hat{u}(\lambda)$
est bien d\'efinie et on a, en posant $h=r^2 g$:
$$
   \forall\,\lambda,\,\Re\lambda=-\varepsilon,\quad
   \rrlb\,\hat{u}(\lambda) = \hat h (\lambda).
$$
 
On rappelle que 
$$
\begin{array}{ll}
   \hat E = H^1(G)\cap\,\bigl\{ W \in H^2(G_-)\times  H^2(G_+)\,,\,
      & a_+\,W_+'(2\pi) - a_-\,W_-'(0) = 0,\\[1mm] 
      & a_+\,W_+'(\frac\pi2) - a_-\,W_-'(\frac\pi2) = 0 \;\bigr\}
   \end{array}
$$
et que $\rrlb$ est inversible de $\hat E$ dans $L^2(G)$ si et seulement si la
condition \Ref{3E1} est v\'erifi\'ee. On voit ainsi que l'on peut choisir
$\varepsilon>0$ assez petit pour que $\rrlb$ soit inversible pour tout $\lambda$,
$\Re\lambda\in[-\varepsilon,0)$.

Gr\^ace au Th\'eor\`eme\REF{final3}, on va d\'emontrer une estimation sur l'inverse
de \rlb:
 \begin{lem}
  \label{new69} 
On suppose que $a_+  a_- \ne -1$.
Soit $\xi\in\R$. Alors il existe $\eta_0$ et $C$ dans $\R_+$ telles que 
$$
   \forall\,V\in \hat E,\,
   \forall\,\lambda = \xi+i\eta\,\mbox{ avec }\,|\eta| \geqslant \eta_0,\quad 
  \Norm{V}{H^2(G_\pm,\,|\eta|)} \leqslant C\,\Norm{\rrlb\,V}{\dg}.
$$
 \end{lem}

\begin{Proof}
On va l'obtenir par la m\'ethode d'addition d'une variable en partant de l'estimation
a priori obtenue dans le Th\'eor\`eme\REF{final3}: pour tout $v\in H_0^1(\Omega)$,
solution du probl\`eme \Ref{pb} avec second membre dans $L^2(\Omega)$:
\begin{equation}
\label{coer4}
   \Norm{v_\pm}{H^2(\Gamma_{1,\pm})} \leqslant 
   C\,\bigl( \Norm{a_+\Delta v_+}{L^2(\Gamma_{2,+})} +
   \Norm{a_-\Delta v_-}{L^2(\Gamma_{2,-})} + \Norm{v}{H^1(\Gamma_2)} \bigr) 
\end{equation}
avec $\Gamma_1$ et $\Gamma_2$ deux couronnes embo\^{\i}t\'ees (par exemple $\Gamma_j$
est l'ensemble des points tels que $r\in(2^{-j},2^j)$). On prend $v=\chi r^\lambda V$,
avec $\chi$ \`a support dans $\Gamma_2$ et valant $1$ sur $\Gamma_1$ et par un calcul
en coordonn\'ees polaires on obtient:
$$
   \Norm{a_\pm\Delta v_\pm}{L^2(\Gamma_{2,\pm})}  \leqslant C\,
   \bigl( \Norm{\rrlb\,V_\pm}{\dg} + \Norm{V}{H^1(G,\,|\eta|)} \bigr),
$$ 
et
$$
   \Norm{v_\pm}{H^2(\Gamma_{1,\pm})}  \geqslant c\,
   \Norm{V_\pm}{H^2(G_\pm,\,|\eta|)} - M\,\Norm{V}{H^1(G,\,|\eta|)},
$$ 
les int\'egrales en $r$ disparaissant par choix du support de $\chi$ et le terme
$\Norm{V}{H^1(G,\,|\eta|)}$ prove\-nant des d\'eriv\'ees de $\chi$. En utilisant
(\ref{coer4}) et les deux in\'egalit\'es pr\'ec\'edentes, on obtient:
$$
   \Norm{V_\pm}{H^2(G_\pm,\,|\eta|)} \le d_1\,\Norm{\rrlb\,V}{\dg} +
   d_2\,\Norm{V}{H^1(G,\,|\eta|)}.
$$
Or 
$$
   \Norm{V}{H^1(G,\,|\eta|)}\leqslant \frac1{|\eta|+1} \bigl(
   \Norm{V_-}{H^2(G_-,\,|\eta|)} + \Norm{V_+}{H^2(G_+,\,|\eta|)} \bigr).
$$
Ainsi pour $|\eta|$ assez grand, $d_2\Norm{V}{H^1(G,\,|\eta|)}$ se majore par la
moiti\'e du membre de gauche de l'avant-derni\`ere estimation, ce qui permet de
prouver le lemme.
\end{Proof}

Nous pouvons en d\'eduire le
\begin{theo}
\label{3T1}
On suppose que $a_+a_-\neq -1$.
Soit $\gamma\in \R$ tel que la condition \Ref{3E1} soit satisfaite sur la droite
$\Re\lambda=1-\gamma$. Alors le probl\`eme de transmission $w \mapsto p$
\begin{equation}
\label{pb5}
   \left\{\begin{array}{ll}
   a_\pm \Delta w_\pm = p_\pm & \mbox{dans }\; \Gamma_\pm,\\[1mm]
   w_+ -  w_- = 0 & \mbox{sur }\; 
   \partial\Gamma_+ \cap \partial\Gamma_-,\\[1mm]
   a_+ \partial_n w_+ - a_- \partial_n w_- = 0 & \mbox{sur }\; 
   \partial\Gamma_+ \cap \partial\Gamma_-,
\end{array}\right.
\end{equation}
induit un isomorphisme de 
$$
   E_\gamma := K^1_{\gamma-1}(\Gamma) \cap 
   \bigl( K_\gamma^2(\Gamma_-)\times K^2_\gamma(\Gamma_+) \bigr)
   \longrightarrow K^0_{\gamma}(\Gamma).
$$
Son inverse $p\mapsto w$ est donn\'e par la formule, {\it cf} \Ref{Mel.inv}
\begin{equation}
\label{Mel.invb}
  \tilde w(r,\theta) =
  \frac{1}{2\pi}\int_{\R}\,r^{\xi+i\eta} \mathcal{L}(\xi+i\eta)^{-1}\hat
   q(\xi+i\eta)\,\ud \eta, \quad\mbox{o\`u $\xi = 1-\gamma$ et $q = r^2p$}.
\end{equation}
\end{theo}

\begin{Proof}
Si $w\in E_\gamma$ satisfait le probl\`eme de transmission \Ref{pb5} avec $p=0$, alors
$$
   \forall\lambda\in\C,\; \Re\lambda = 1-\gamma,\quad
   \rrlb\,\hat{w}(\lambda) = 0.
$$
L'unicit\'e provient alors de l'injectivit\'e de \rlb sur la droite $\Re\lambda =
1-\gamma$. 

Appliquons le lemme\REF{new69} \`a la fonction $V(\lambda) :=
\rrlb^{-1}\hat{q}(\lambda)$ avec
$\Re\lambda=\xi=1-\gamma$:
$$
   \Norm{V(\lambda)}{H^2(G_\pm,\,|\Im\lambda|)} \leqslant
   C\Norm{\rrlb\,V(\lambda)}{\dg},
$$ 
pour $|\Im\lambda|>\eta_0$. Mais pour $|\Im\lambda|\leqslant\eta_0$, l'inversibilit\'e de \rlb
et la continuit\'e de l'inverse par rapport \`a $\lambda$ permet aussi d'avoir
l'estimation ci-dessus. D'o\`u:
$$
   \int_\R  \Norm{V(\xi+i\eta)}{H^2(G_\pm,\,|\eta|)}^2 \,\ud\eta \leqslant C\int_\R
  \Norm{\hat q(\xi+i\eta)}{\dg}^2 \,\ud\eta.
$$ 
Le majorant est fini, car
$(\eta,\theta)\mapsto \hat{q}(\xi+i\eta)\in L^2(\R\times G)$ par le lemme\REF{5.4}.
Donc 
$$
   (\eta,\theta) \mapsto V(\xi+i\eta,\theta) \in L^2(\R,H^2(G_\pm)),
$$ 
et on conclut alors, toujours par le lemme\REF{5.4}.
\end{Proof}

Nous avons imm\'ediatement les deux cons\'equences suivantes pour notre probl\`eme
de transmission initial \Ref{pb4}:  
\begin{cor}
  \label{reg.Melu}
 Pour $\varepsilon>0$ tel que la condition \Ref{3E1} soit satisfaite sur la droite
$\Re\lambda=-\varepsilon$, la solution $u$ du probl\`eme \Ref{pb4} est
telle que $u$ appartient \`a $K_{1+\varepsilon}^2(\Gamma_\pm)$.
\end{cor}

\begin{Proof}
Il suffit d'appliquer le Th\'eor\`eme\REF{3T1} pour $\gamma = 1+\varepsilon$, car
comme le second membre $g$ est $L^2(\Gamma)$ \`a
support compact, il appartient \`a $K^0_{1+\varepsilon}$. Comme $u$ est d\'ej\`a dans
$K^1_{\varepsilon}(\Gamma)$, elle co\"{\i}ncide avec la solution
$w\in E_{1+\varepsilon}$ du probl\`eme \Ref{pb5} avec second membre $p=g$.
\end{Proof}

Et gr\^ace au Lemme\REF{3L3}:
\begin{cor}
  \label{reg.Melu0}
Il existe une unique solution $u_0\in E_0$ au probl\`eme
de transmission \Ref{pb4}.
\end{cor}

\begin{Proof}
Vu que la condition \Ref{3E1} est satisfaite sur la droite $\Re\lambda=1$, il suffit
cette fois-ci d'appliquer le Th\'eor\`eme\REF{3T1} pour $\gamma = 0$.
\end{Proof}

Il reste \`a comparer $u$ et $u_0$. Avec $h = r^2g$, en prenant
$\lambda = -\varepsilon - i\eta$ comme param\'etrisation de la droite
$\Re\lambda = -\varepsilon$ et  $\lambda = 1+i\eta$ comme param\'etrisation de la
droite $\Re\lambda=1$, (\ref{Mel.invb}) donne:
\begin{equation}
\label{Mel.invc}
\begin{array}{l}
  \displaystyle\tilde u(r,\theta) = -
  \frac{1}{2i\pi} \downarrow \int_{\Re\lambda=-\varepsilon}\,r^{\lambda}\rrlb^{-1}\hat
   h(\lambda)\,\ud \lambda,\\[4mm] \displaystyle\mbox{}\hskip 2cm\mbox{et}\quad
  \tilde u_0(r,\theta) =
  \frac{1}{2i\pi} \uparrow \int_{\Re\lambda=1}\,r^{\lambda}\rrlb^{-1}\hat
   h(\lambda)\,\ud \lambda,
\end{array}
\end{equation}
avec $h = r^2g$. On va maintenant expliquer pourquoi, comme il est classique en
th\'eorie des probl\`emes \`a coins, voir \cite{Kondratev}, la diff\'erence $u-u_0$
s'exprime comme une somme finie de r\'esidus d'une fonction m\'eromorphe.

\textbf{1.} Le symbole $\mathcal{L}$ est polynomial en $\lambda$ \`a valeurs dans
l'espace des op\'erateurs born\'es de $\hat E$ dans  $L^2(G)$ (voir d\'ebut de la
section\REF{3S3}), donc analytique. Son inverse existe sauf en un ensemble
d\'enombrable de points en vertu du Lemme\REF{3L2}. On d\'eduit du th\'eor\`eme
analytique de Fredholm (et en remarquant que l'injection $\hat E \hookrightarrow
L^2(G)$ est compacte) le r\'esultat suivant:

\begin{theo}  
   \label{mero}
L'inverse du symbole Mellin $\lambda \mapsto \rrlb^{-1}$ est une fonction m\'eromorphe
sur $\C$ \`a valeurs dans l'espace $\mathbb{L}(L^2(G),\hat E)$ des applications
continues de $L^2(G)$ dans $\hat E$. De plus ses parties polaires sont de rang fini.
L'ensemble de ses p\^oles est not\'e $\sigll$, qui est le spectre de $\mathcal{L}$.
\end{theo}   

\textbf{2.} La transform\'ee de Mellin du second membre est analytique:

\begin{lem}
  \label{gholo} 
Avec $h = r^2g$, o\`u $g$ est d\'efinie au d\'ebut de la section\REF{3S3}, la
fonction $\hat h$ est holomorphe \`a valeurs dans $L^2(G)$ \`a l'int\'erieur de la
bande $-\varepsilon<\Re\lambda<1$ pour tout $\varepsilon>0$.
\end{lem}

\begin{Proof}
Il suffit de montrer que 
$$
   \int_0^\infty\, \Norm{r^{-\lambda}\hat h(r)}{L^2(G)} \frac{\ud r}{r} < \infty.
$$
En coupant l'int\'egrale en deux int\'egrales, de bornes $0$ et $1$, resp.\ $1$ et
$+\infty$, on majore l'int\'egrale par $c(\Norm{g}{\kkk} + \Norm{g}{\kkkk})$, o\`u
les deux normes sont born\'ees du fait que $g$ est dans $L^2(\Gamma)$ \`a support
compact.
\end{Proof}

On choisit $\varepsilon$ telque $\sigll\cap\{-\varepsilon<\Re\lambda<0\}$ soit vide. 
Comme $\rrlb^{-1}$ est m\'eromorphe et $\hat{h}(\lambda)$ holomorphe, pour tout
$r>0$ fix\'e, la fonction $r^{\lambda}\rrlb^{-1}\,\hat{h}(\lambda)$ est m\'eromorphe
\`a valeurs dans $\hat E$ sur la bande $-\varepsilon<\Re\lambda<1$. L'ensemble $\sigl$
d\'efini par $\sigll\cap\{-\varepsilon<\Re\lambda<1\}$ est fini d'apr\`es le
Lemme\REF{new69} et le Th\'eor\`emes\REF{mero}. De plus par choix de $\varepsilon$
$$
   \sigl = \sigll\cap\{0\leqslant \Re\lambda<1\}.
$$
Soit $\mathcal{C}$ un contour simple entourant $\sigl$ et contenu dans 
$-\varepsilon<\Re\lambda<1$. Par le th\'eor\`eme de Cauchy:
\begin{equation}
\label{Cau1}
   \frac{1}{2i\pi}\,\int_\mathcal{C} r^\lambda \rrlb^{-1}\hat h(\lambda)\,\ud \lambda
   = \sum_{\lambda_0\,\in\, \sigl} \;\Res_{\lambda\,=\,\lambda_0}\,
   r^\lambda \rrlb^{-1}\hat h(\lambda).
\end{equation}
Le Lemme\REF{new69} et des estimations sur les transform\'ees de Mellin permettant de
montrer\\[-22pt] 
\begin{itemize}
  \item qu'on peut pousser le contour $\mathcal{C}$ jusqu'au bord de la bande
        $-\varepsilon<\Re\lambda<1$,\\[-20pt]
  \item que les c\^ot\'es horizontaux du rectangle infini ne contribuent pas \`a
        l'int\'egrale,\\[-20pt]
\end{itemize}
on d\'eduit de \Ref{Mel.invc} la formule:
\begin{equation}
\label{Cau3}
   u_0 - u = \sum_{\lambda_0\,\in\, \sigl} \;
   \Res_{\lambda\,=\,\lambda_0}\, r^\lambda \rrlb^{-1}\hat h(\lambda). 
\end{equation}

\subsection{D\'ecomposition de la solution}
Il reste \`a exploiter la formule ci-dessus: expliciter les r\'esidus autant que
faire ce peut et retourner au probl\`eme initial sur $\Omega$.

On \'etudie d'abord un peu plus en d\'etail l'ensemble $\sigl$. Il s'agit de trouver
les racines $\lambda$ de l'\'equation
\begin{equation}
\label{3E1b}
   (\mu+1)\,\sin\lambda\pi \;=\; \pm\, (\mu-1)\,\sin\lambda(\pi-\omega),
\end{equation}
contenues dans la bande $\Re\lambda\in[0,1]$. La valeur particuli\`ere de l'angle
$\omega = \pi/2$ permet des calculs explicites ais\'es. Posons $t = \lambda\pi/2$.
Alors l'\'equation \Ref{3E1b} devient
$$
   2(\mu+1)\,\sin t\cos t \;=\; \pm\, (\mu-1)\,\sin t.
$$
D'o\`u les racines $t = k\pi$, correspondant \`a $\lambda=2k$, $k\in\Z$, et
l'\'equation r\'esiduelle
$$
   \cos t = \pm\rho,\quad\mbox{avec}\quad
   \rho = \frac{|\mu-1|}{2|\mu+1|}\,.
$$
Cette \'equation a les racines $t = \arccos\rho + 2k\pi$ et $t = -\arccos\rho +
(2k-1)\pi$ lorsque $\rho\leqslant 1$, et les racines $t = 2k\pi \,\pm\, i\argch\rho$
lorsque $\rho>1$. D'autre part $\rho>1$ si et seulement si $-3<\mu<-1/3$. Ceci permet
imm\'ediatement d'\'etablir

\begin{lem}
\label{3L4}
On suppose $\mu<0$ et $\mu\neq-1$.\\[-22pt] 
\begin{itemize}
\item Si $\mu\geqslant-\frac13$ ou $\mu\leqslant -3$, 
\begin{equation}
\label{3E4a}
   \sigl = \{0, \lambda_1(\mu)\} 
\end{equation}
o\`u la fonction $\lambda_1$ est croissante sur $[-\frac13,0]$, avec $\lambda_1(-\frac13)
= 0$ et $\lambda_1(0) = \frac23$; enfin $\lambda_1(\frac1\mu) =
\lambda_1(\mu)$.\\[-20pt] 
\item Si $-3<\mu<-\frac13$, 
\begin{equation}
\label{3E4b}
   \sigl = \{0, \pm\, i\eta(\mu)\} 
\end{equation}
o\`u la fonction $\eta(\mu)$ est positive, d\'ecroissante sur $(-1,-\frac13]$, tend
vers l'infini quand $\mu\to-1$, $\eta(-\frac13)=0$ ; enfin $\eta(\frac1\mu) =
\eta(\mu)$.
\end{itemize}
\end{lem}

  Pour $\lambda_0\in \sigll$, introduisons l'espace
$$
   Z^{\lambda_0} = \bigl\{ \Res_{\lambda\,=\,\lambda_0}\, r^\lambda
   \rrlb^{-1}G(\lambda) , 
   G(\lambda)\mbox{ holomorphe en }\lambda_0\bigr\} .
$$
Un d\'eveloppement en s\'erie de $r^\lambda$ au voisinage de $\lambda_0$ donne
imm\'ediatement que les \'el\'e\-ments de $Z^{\lambda_0}$ ont la forme d'une somme
finie de termes en $r^{\lambda_0}\log^q \!r \,\varphi(\theta)$. La formule \Ref{Cau3}
se r\'e\'ecrit alors
\begin{equation}
\label{Cau3b}
   u_0 - u = \sum_{\lambda_0\,\in\, \sigl} \;
   \Res_{\lambda\,=\,\lambda_0}\, z_{\lambda_0} \quad\mbox{o\`u}\quad
   z_{\lambda_0} \in Z^{\lambda_0}. 
\end{equation}
Tronquons par $\chi$ la d\'ecomposition ci-dessus. On a donc $\chi u_0 \in
H^2(\Omega_\pm)$ et, comme $\chi u$ est par hypoth\`ese dans $H^1(\Omega)$, la somme
des $\chi z_{\lambda_0}$ doit aussi appartenir \`a $H^1(\Omega)$. Vu leur structure
en $r^{\lambda_0}$, cela exige que chacun d'entre eux soit dans $H^1(\Omega)$. Pour
$\lambda_0$ de partie r\'eelle $>0$, c'est toujours vrai. Par contre, si
$\Re\lambda_0=0$ sans que $\lambda_0$ soit nul, cela implique que $z_{\lambda_0}=0$.
Enfin, si $\lambda_0=0$, cela implique que $z_{\lambda_0}$ soit une constante, donc
une fonction r\'eguli\`ere. C'est pourquoi l'on pose finalement
$$
   u_\reg = \chi (u_0 - z_0) \in H^2(\Omega_\pm),
$$
et que \Ref{Cau3b} devient
\begin{equation}
\label{Cau3c}
   u_\reg - \chi u = \sum_{\lambda_0\in \sigll,\; 0<\Re\lambda<1}
   \;\Res_{\lambda\,=\,\lambda_0}\, z_{\lambda_0} \quad\mbox{o\`u}\quad
   z_{\lambda_0} \in Z^{\lambda_0}. 
\end{equation}
Ceci, avec les pr\'ecisions fournies par le Lemme\REF{3L4} permet finalement
d'obtenir:

\begin{theo}
 \label{fin}
Soit $u\in H^1_0(\Omega)$ solution du probl\`eme \Ref{pb} avec second membre $f\in
L^2(\Omega)$. Soit $\mathcal{O}$ un coin de transmission (un coin de $\Omega_-$) et
$\chi$ une fonction de troncature qui vaut $1$ au voisinage de $\mathcal{O}$. On pose
$\mu = a_+/a_-$ et on suppose que $\mu<0$ et que $\mu\neq-1$. Alors:
\\[-22pt] 
\begin{itemize}
\item Si $\mu>-\frac13$ ou $\mu< -3$, 
\begin{equation}
\label{3E8a}
   \chi u = u_\reg + c_1 \chi r^{\lambda_1}\varphi_1(\theta)
\end{equation}
o\`u la fonction $u_\reg \in H^1(\Omega) \cap \bigl(H^2(\Omega_-)\times
H^2(\Omega_+)\bigr)$, l'exposant $\lambda_1$ est d\'efini dans le
Lemme\REF{3L4} et la fonction $\varphi_1$ est dans $H^1(G)$ avec
$\varphi_{1,\pm}\in \mathcal{D}(\overline G_\pm)$.\\[-20pt] 
\item Si $-3\leqslant\mu\leqslant-\frac13$, 
\begin{equation}
\label{3E8b}
   \chi u_\pm \in H^2(\Omega_\pm).
\end{equation}
\end{itemize}
\end{theo}


\section{Conclusion et perspectives}

Apr\`es les th\'eor\`emes\REF{final3} et\REF{fin}, la question qui reste
pendante est celle de l'existence d'une solution au probl\`eme \Ref{pb}.

Par le th\'eor\`eme de Lax-Milgram, apr\`es multiplication par un nombre
complexe bien choisi, on montre facilement que le probl\`eme \Ref{pb} a une
unique solution si $a_+/a_-$ n'appar\-tient pas \`a $\R_-$. Cela permet de
d\'eduire que
\medskip

\begin{center}
\noindent
\textit{Si l'op\'erateur $A$ du probl\`eme \Ref{pb} est \`a indice\footnote{Un
op\'erateur est dit \`a indice si son noyau et son conoyau sont de dimension
finie, l'indice est alors la dimension du noyau moins celle du conoyau.
Rappelons en particulier qu'un op\'erateur \`a indice est \`a \textit{image
ferm\'ee.} De plus l'indice est une fonction continue dans l'ensemble des
op\'erateurs \`a indice pour la topologie des op\'erateurs born\'es.}
de $H^1_0(\Omega)$ dans $H^{-1}(\Omega)$,\\ alors son indice est nul.}
\end{center}

\medskip
Il reste donc \`a savoir si l'op\'erateur $A$ est \`a indice. La m\'ethode
d'investigation envisag\'ee est la construction d'une param\'etrix (i.e.\ un
inverse modulo un op\'erateur compact). Une telle construction s'effectue par
localisation, la seule difficult\'e restante \'etant le voisinage de chaque
coin de transmission int\'erieur, \`a cause de la pr\'esence permanente de
p\^oles du symbole $\rrlb^{-1}$ sur la droite $\Re\lambda=0$.

La conjecture est la suivante :
\begin{itemize}
\item Si $\mu>-\frac13$ ou $\mu< -3$, l'op\'erateur $A$ est \`a indice.\\[-20pt] 
\item Si $-3\leqslant\mu\leqslant-\frac13$, l'op\'erateur $A$ n'est pas \`a
indice (car \`a image non ferm\'ee).
\end{itemize}

Ceci repose sur une analyse fine de $\rrlb^{-1}$ au voisinage de ses p\^oles
situ\'es sur la droite $\Re\lambda=0$. Dans le premier cas, le seul p\^ole dans
cette situation est $\lambda=0$: il est double, donnant un espace $Z^0$
engendr\'e par $1$ est $\log r$; $1$ est une fonction r\'eguli\`ere et le
coefficient devant $\log r$ dans l'asymptotique est une forme lin\'eaire
continue sur l'espace des seconds membres (la dualit\'e contre la fonction
constante $1$), contrairement \`a la situation du second cas.



\end{document}